\documentclass{article}
\usepackage{amssymb,amsmath,amsfonts}
\newtheorem{theorem}{Theorem}
\newtheorem{example}{Example}
\newtheorem{lemma}{Lemma}
\newtheorem{remark}{Remark}
\newtheorem{corollary}{Corollary}
\newtheorem{proposition}{Proposition}

\title{On possible growths of Toeplitz languages}
\author{J. Cassaigne\footnote{Institut de Math\'ematiques de Luminy 
13288 Marseille Cedex 9, France; \tt{cassaigne@iml.univ.mrs.fr}}, A. Frid\thanks{Supported in part by RFBR grants 09--01--00244 and 10--01--00424}\footnote{Sobolev Institute of Mathematics SB RAS, 4 Koptyug av. 630090 Novosibirsk; \tt{anna.e.frid@gmail.com}}, and F. Petrov\thanks{Supported in part by RFBR grant 08--01--00379 and the grant of President of Russia for leading scientific schools no. NSh--2460.2008.1 }\footnote{St-Petersburg division of Steklov Mathematical Institute, 27 Fontanka emb., 191023, St.-Petersburg; {\tt fedyapetrov@gmail.com} }}
\date{Submitted March 10 2010}

\begin{document}
\maketitle
\begin{abstract}
We consider a new family of factorial languages whose subword complexity grows as $\Theta(n^{\alpha})$, where $\alpha$ is the only positive root of some transcendent equation. The asymptotic growth of the complexity function is studied with the use of analytical methods and in particular with a corollary of the Wiener-Pitt theorem. The factorial languages considered are also languages of arithmetic factors of infinite words; so, we describe a new family of infinite words with an unusual growth of arithmetical complexity.

\end{abstract}

\section{Introduction}
This study started as an attempt to construct an infinite word with a non-standard arithmetical complexity growth. Estimation of this complexity with the use a recurrent relation required non-trivial and non-elementary analytical methods, namely, the Wiener Tauberian theory. Thus, apart from the main combinatorial result (Theorem 1), its analytic derivation is also interesting. 

The considered family of words is built by a rather simple construction, and their arithmetical complexity grows as $\Theta(n^{\alpha})$, where $\alpha$ is a root of some transcendental equation. This is proved by elementary techniques (Lemma 4). Then to give the precise asymptotics we use a corollary of the Wiener-Pitt theorem (Theorem 3).

\section{Definitions, survey and results}
We consider finite and infinite words on the alphabet $\Sigma$ of cardinality $d\geq 2$. A set of finite words on it is called a {\it language}; in particular, we can speak on the language of factors Fac$(w)$ of an infinite word $w$. Such a language is always {\it factorial}, which means that it is closed under taking factors; here {\it factor}, or {\it subword} $v$ of a word $u$ is a word such that $u=svt$ for some $s$ and~$t$. In general we can speak on the {\it factorial closure} of a language which is the set of factors of all its elements.

An {\it arithmetical factor} of a word $w_1\dots w_n \dots$ with $w_i \in \Sigma$ is a finite  word of the form $w_kw_{k+r}w_{k+2r}\dots w_{k+mr}$ with $r>0$. The set of all arithmetical factors of a given word is its {\it arithmetical closure}. It is clear that the arithmetical closure of a word or of a language is a factorial language.

The number of words of length $n$ in a factorial language $F$ is called its {\it subword complexity} and is denoted by $p_F(n)$; if this factorial language is the language of factors of an infinite word $w$, then its subword complexity is denoted by $p_w(n)$. The subword complexity of the arithmetical closure of $w$ is called its {\it arithmetical complexity} and denoted by  $a_w(n)$.

It is clear that both functions of an infinite word do not decrease and grow not faster than $d^n$, where $d$ is the cardinality of the alphabet. A survey on the subword complexity can be found in [1]; a more recent result to be mentioned is the construction of a word whose subword complexity grows faster than any polynomial and slower than any exponent [2].

The function of arithmetical complexity introduced in [3] for a non-periodic word can grow both exponentially [4] and linearly; moreover, a characterization of uniformly recurrent words of linear arithmetical complexity is known [5]. Then, there are some examples of words whose arithmetical complexity grows not linearly but subpolynomially [6], or faster than any polynomial but slower than  $c^n$ for any $c>1$ [7].

The goal of this paper is the construction of a family of infinite words whose arithmetical complexity grows in a non-standard way, namely, as $n^\alpha$, where $\alpha$ is a root of some transcendental equation.

To state the main result of the paper, let us fix a finite set of prime numbers $Q=\{q_1,\ldots,q_k\}$ of cardinality $k\geq 2$ and a number $d$ which is the size of the alphabet considered. Consider the equation
\begin{equation}\label{e1}
\frac{d-1}{d}=\prod\limits_{i=1}^{k} ( 1-q_i^{-x} )
\end{equation}
and note that for each $Q$ and $d$ it has the unique positive root $x=\alpha=\alpha(Q,d)$. Indeed, the function $\Pi(x)=\prod\limits_{i=1}^{k} \bigl( 1-q_i^{-x} \bigr)$ is growing with $x \geq 0$, and we have $\Pi(0)=0$ and $\Pi(x) \to 1$ with $x \to \infty$.

Note that for each fixed $d$ the root $\alpha(Q,d)$ of the equation (1) never exceeds the root $\alpha=\alpha(d)$ of the equation
$$
\frac{d-1}{d}=\prod\limits_{q \mbox{~is~prime}} (1-q^{-x})=
\frac{1}{\zeta(x)},
$$
where $\zeta$ is the Riemann zeta-function. In particular, if $d=2$, then for any set $Q$ the value of $\alpha$ is not greater than  the root $\alpha=1.7286\ldots$ of the equation $\zeta(x)=2$.

As it follows from the Besikovitch theorem [8], the root $\alpha$ of (1) can be either integer (like in the case of $d=3$ and $Q=\{2,3\}$, when $\alpha=2$), or irrational. The Schanuel conjecture [9] would imply that each irrational root of this equation is transcendental; however, to the best of our knowledge, at the moment this fact is not proven.

The main combinatorial result of this paper is the following
\begin{theorem}
Let $Q=\{q_1,\ldots,q_k\}$ be a finite set of prime numbers, $k\leq 2$. Then for each $d>1$ there exists an infinite word $w=w(Q,d)$
on the $d$-letter alphabet such that its arithmetical complexity $a_w(n)=p_{Q,d}(n)$ satisfies
$$
\frac{p_{Q,d}(n)}{n^{\alpha(Q,d) +1}} \to C \quad \mbox{~with~} n \to \infty,
$$
where $C$ is a positive constant.
\end{theorem}
In what follows we prefer to consider $p_{Q,d}(n)$ as the subword complexity function of the language $F=F(Q,d)$ equal to the arithmetical closure of $w$. The language $F$ and the infinite word $w$ are constructed in the next section. Then in Section 4 the classical method of special words is used to derive a recurrent formula for the first differences $s_{Q,d}(n)=p_{Q,d}(n+1)-p_{Q,d}(n)$ of the function $p_{Q,d}(n)$. The properties of the function $s_{Q,d}(n)$ are studied in Section 5 by combinatorial techniques. Then in Section 6 we prove a property of $\Pi(x)$ as a complex funcion. That property allows to apply to it Theorem 3 from Section 7, which is a corollary of a Tauberian theorem. The application is made in Section 8 which completes the proof of Theorem 1.

\section{The language $F(Q,d)$}
Let $\Sigma$ be a finite alphabet of cardinality $d$. Recall that a {\it morphism} is a mapping $\varphi: \Sigma^*
\to \Sigma^*$ which for all words $x,y \in \Sigma^*$ satisfies the equality $\varphi(xy)=\varphi(x)\varphi(y)$; here and below two words written consecutively mean their concatenation.

Now for all $a \in \Sigma$ and $m \in \mathbb N$ define the morphism
 $\varphi_{a,m}$ by the equalities
$$
\varphi_{a,m}(b)=a^{m-1}b
$$
for all $b \in \Sigma$ (in particular, for $b=a$). Note that this morphism can be interpreted also as a Toeplitz transform, since the word $\varphi_{a,m}(x)$ is obtained by substituting consequtive symbols of $x$ to the ``holes'', denoted by diamonds, of the partial word $a^{m-1}\diamond a^{m-1}\diamond \dots$ [10,\,11].

Let us fix a finite set $Q$ of prime numbers, $Q=\{q_1,\ldots,q_k\}$. Define the language $L=L(\Sigma,Q)$ as the closure of
the alphabet $\Sigma$ with respect to all the morphisms $\varphi_{a,q}$, where
 $a \in \Sigma$ and $q \in Q$:
$$
L(\Sigma,Q)=\{\varphi_{c_m,r_m}(\varphi_{c_{m-1},r_{m-1}}(\dots(
\varphi_{c_1,r_1}(c_0))\dots ))\mid c_i \in \Sigma,\ r_i \in Q\}.
$$

In what follows we study the factorial closure $F=F(\Sigma,Q)=$Fac$(L(\Sigma,Q))$ of the language $L(\Sigma,Q)$. This factorial language obviously is not equal to the language of factors of any infinite word, but is equal to the language of arithmetical factors of the word $w$ defined as the limit  
$$
w=\lim_{n\to\infty}
\varphi_{c_1,r_1}(\varphi_{c_2,r_2}
(\dots(\varphi_{c_n,r_n}(c_{n+1}))\dots )),
$$
where the sequence of pairs  $\{(c_i,r_i)\}_{i=1}^{\infty}$ for all $j$ contains all the $(dk)^j$ factors of length $j$, defined as the sequences from  $(\Sigma \times Q)^j$. It is not difficult to see that the limit $w$ exists, and its arithmetical closure is equal to $F$; in particular, each arithmetical subsequence of difference  $r_1\dots r_m$ either consists of equal symbols or contains all the factors of the form 
$\varphi_{c_{m+1},r_{m+1}}(\varphi_{c_{m+2},r_{m+2}}
(\dots(\varphi_{c_{m+n},r_{m+n}}(c_{m+n+1}))\dots ))$.

The complete proof of the fact that the arithmetical closure of $w$ is $F$ is analogous to the proofs from [6], where the case of $Q=\{p\}$ (that is, $|Q|=1$) was considered: the order of growth of the arithmetical complexity in that case is equal to $\Theta(n^{1+\log_p d})$, but the limit of the ratio
of $p_{\{p\},d}(n)$ to $n^{1+\log_p d}$ does not exist. In the same paper, the arithmetical complexity of other words generated by Toeplitz constructions with the unique length of the morphism was considered, so that we can focus on the case of  $|Q|\geq 2$.
\begin{example}{\rm
Let $\Sigma=\{a,b\}$ and $Q=\{3,5\}$; consider the language $F=F(\Sigma, Q)$. The word 
$u=abaabaabaabaaaaabaaba$ belongs to $F$ since it is a factor of the word $\varphi_{a,3}(bbbbabbb)$; then, $bbbbabbb$ is a factor of $\varphi_{b,5}(aa)$, and $aa$ is a factor of $\varphi_{a,3}(a)$. Thus, 
$u \in $Fac$(\varphi_{a,3}(\varphi_{b,5}(\varphi_{a,3}(a))))\subset F$.

Note that in fact, we can uniquely reconstruct from $u$ only the last morphism applied. In fact, the word $bbbbabbb$ is also a factor of $\varphi_{b,3}(bab)$, and $bab$ is a factor of both $\varphi_{b,3}(a)$ and 
$\varphi_{b,5}(a)$.
}
\end{example}
\begin{remark}{\rm 
For all $n,m>0$ and $a \in \Sigma$ the morphisms $\varphi_{a,n}$ and
$\varphi_{a,m}$ commute:
$\varphi_{a,n}\varphi_{a,m}=\varphi_{a,nm}=\varphi_{a,m}\varphi_{a,n}$.
}
\end{remark}

The remaining part of the work will be related to finding the asymptotics of the function $p_F(n)=p_{d,Q}(n)$ ($p(n)$ for short). Most arguments and computations will concern not $p(n)$ itself but its dirst differences 
$$
s(n)=p(n+1)-p(n).
$$

\section{Recurrent formula on first differences}
It is well known (see e.~g.[12]) that the first differences $s(n)$ of the complexity function can be expressed in terms of specialwords of the factorial language $F$.

Consider a word $u \in F$ and denote by $L(u)$ the set of symbols 
$a \in\Sigma$ such that $au \in F$. The cardinality of $L(u)$ is called the (left)
{\it speciality degree} of $u$ and denoted by $l(u)$. A word
$u$ with $l(u)\neq 1$ is called {\it special} in $F$; the set of all special words of length $n$ in $F$ is denoted by  $S(n)$. The following formula expressing the complexity of $F$ in terms of its special words is well-known: 
\begin{equation}
s(n)=\sum\limits_{u\in S(n)}(l(u)-1).
\end{equation}
So, to find the expression for the first differences of $F$, we should first study the set of its special words.

\begin{lemma}
Consider a special word $u\in S(n)$ of $F$ of length at least three. Then one of the following is true:

$(1)$~either $u=a^n$ for some symbol $a \in \Sigma$,

$(2)$~or $u=\varphi_{a,N}(v)a^{N'}$ for some
$a \in \Sigma$, $N \in {\mathbb N}$, where all the prime divisors of $N$
belong to $Q$, $0\leq N'<N$ and $v \in F$; moreover, $v$
contains the symbol $b \neq a$ and $L(v)=L(u)$.
\end{lemma}
{\sc Proof.} Consider first the situation when $2 \notin Q$.
In this case each word $u$ from $F$ of length at least 3 either contains two consecutive equal symbols, or is equal to 
$u=xyx$ for some $x \neq y$, $x,y \in \Sigma$.
In the latter case $u$ occurs in $F$ only as a factor of words of the form  $\varphi_{x,q}(v)$ for $q \in Q$ and $v \in F$, and is not special since the only letter by which $u$ can be extended to the right to an element of $F$ is $x$. 

Thus, in each word of length at least 3 which is special in $F$ there are two consecutive equal symbols. If we denote it by $a$, we see that $u$ occurs in $F$ only as a factor of words of the form  
$\varphi_{a,q}(v)$, where $q \in Q$ and $v \in F$.

The situation when $u=a^n$ corresponds exactly to the case (1) of the lemma. Note that
$L(a^n)=\Sigma$: indeed, each extension $xa^n$ of $a^n$, where $x\in
\Sigma$, is in $F$, since it is a factor of the word 
$\varphi_{a,q}^m(xx)$, and thus of the word
$\varphi_{a,q}^m \varphi_{x,q}(x)$ for all $q \in Q$ and sufficiently large $m$.

Now suppose that $u \neq a^n$. Since $u$ is special, it can be extended to the left with some symbol $b \neq a$. The word $bu$ is a factor of some word $\varphi_{a,q}(v')$, where $q \in Q$ and $v' \in F$. So, all distances between the symbols of  $bu$ not equal to $a$ are divisible by at least one number $q \in Q$. Since all elements of $Q$ are (co)prime, in fact the minimal distance between to consecutive symbols not equal to $a$ is divisible by some product of powers of elements of $Q$: let us denote it by $N(u)=q_1^{l_1}\dots q_k^{l_k}$. Here all $l_i $ are non-negative and at least one of them is positive. Thus, 
$u=\varphi_{a,N(u)}(v)a^{N'}$, and if the word $v$ (for this fixed $N(u)$) is chosen to be the longest possible, then
$0\leq N' <N(u)$. Note that since $N(u)$ is chosen to be maximal, the word $v$ does not contain two consecutive $a$s. 

It remains to prove that $L(v)=L(u)$. The inclusion $L(v) \subseteq L(u)$ is nearly obvious 
since for each $b \in \Sigma$ the fact that $bv \in F$ implies that 
$\varphi_{a,N(u)}(bv)a^{N'} \in F$, and thus that $bu \in F$ (since $a^{N'}$
is a prefix of the $\varphi_{a,N(u)}$-image of any extension of $bv$ to the right).

Let us prove the opposite inclusion $L(u) \subseteq
L(v)$ by the induction on $l=\sum\limits_{i=1}^k l_i$ Suppose first that $l=1$, that is, $N(u)=q \in Q$. Then for each symbol  $x \in \Sigma$, $x\neq a$ the word $xu$ can occur in $F$ only as a suffix of some word $\varphi_{a,q}(xv)a^{N'}$, where $v$
does not contain two consecutive  $a$s. Here by the construction we see that if $x \in L(u)$, then $x \in L(v)$. The same argument works for $x=a$ if $v$ contains at least two symbols not equal to $a$, and the number $q$
can be found from the distance between them. At last, if $v$
contains only one symbol $c$ not equal to $a$, that is, if 
$v$ is a factor of $aca$, then $a \in L(u)$ and $a \in L(v)$. The base of induction is proved.

To prove the insuction step, consider the general case:
$N(u) =q_1^{l_1}\dots q_k^{l_k}$; without loss of generality suppose that 
$l_1,\dots,l_i>0$ and $l_{i+1}=\dots=l_k=0$. So, $u$ can occur in $F$ as a result of application of any of the morphisms  $\varphi_{a,q_j}$,
$j \leq i$,
to the word $\varphi_{a,N(u)/q_j}(v)$, 
and adding the word $a^{N'}$ to the right (which is always possible). But for all  $j\leq i$ by the induction hypothesis we have  $L(\varphi_{a,N(u)/q_j}(v))=L(v)$, that is,
$L(u)=L(v)$, which was to be proved. 

Note that here we have used two facts: first, that the morphisms 
$\varphi_{a,p}$ and  $\varphi_{a,q}$ commute for all $p$ and $q$, and second, that the elements of $Q$ are (mutually co)prime: due to this, the inverse image of $u$ under any morphism $\varphi_{a,q_j}$, $j\leq i$, is really an element of $F$.

For the case of $2 \notin Q$, the lemma is proved. Suppose now that $2 \in Q$.
In this case, as before, if the word contains two consecutive equal letters, then this letter correponds to the last morphism applied, and thus we can apply all the arguments from the previous case. The only difference is that there are long words in $F$ in which there are no two consecutive equal letters. Such a word can be special only if it is of the form
$u=av_1av_2a\dots$
for some letters $a, v_i \in \Sigma$, that is, it can be represented as 
$u=\varphi_{a,2}(v_1v_2\dots)a^m$, where $m\in\{0,1\}$.

If some two of the symbols $v_i$ are distinct, then obviously the last morphism applied correspond not to them but to $a$. So, we can treat the word $u$ as in the previous case. And if all $v_i$ are equal, that is, if $u=ababa\dots$
for some symbol $b \neq a$, then $u$ is the image under $\varphi_{a,2}$ of the word $b^n$ for the respective $n$ (or a prefix of such an image if $m=1$). Thus, $L(u)=L(v)=\Sigma$. Note that at the same time, $u$ is a factor of $\varphi_{2,b}(a^{n'})$ for some $n'$, but it cannot add anything to the set of extensions of $u$ since that set is already maximal. 

\begin{corollary}
For each special word $u \in S(n)$ of length $n \geq 3$ in $F$ there exist a unique letter
$a$ and a unique subser  $Q(u)$ of $Q$ such that for all $q \in Q(u)$ 
{\rm (} and for no other elements of $Q)$
the word $u$ can be represented as $u=\varphi_{a,q}(v_q)a^m$ for the respective $v_q \in F$ and $0\leq m <q$. The length of the word $v_q$ is here equal to 
$\lfloor n/q \rfloor$, and for all $q\in Q(u)$ the equalities hold $L(u)=L(v_q)$.
\end{corollary}
{\sc Proof.} Special words containing at least two distinct symbols are considered in the case (2) of the previous lemma. For the words of the form $a^n$ we have
$L(u)=\Sigma$, $Q(u)=Q$ and $v_q=a^{\lfloor n/q \rfloor}$:
in particular, if the length of $u$ is less then some $q \in Q$, then $v_q$ is the empty word. 

\medskip
Now note that we can apply to each special word  $v \in S(m)$ each morphism $\varphi_{a,q}$ with $a \in\Sigma$ and $q \in Q$, and then add any suffix  $a^r$ for $r<q$. Due to Lemma 1, the obtained word $u=\varphi_{a,q}(v)a^r$ of length $mq+r$ will be special with $L(u)=L(v)$, and we have shown that all the special words of length at least 3 can be obtained that way. Here a word $u$ can have several inverse images under different morphisms (but based on the same letter $a$): $u=\varphi_{a,q_i}(v_i)a^{r_i}=\varphi_{a,q_j}(v_j)a^{r_j}$ for some
$i \neq j$ and respective $v_i,v_j,r_i,r_j$. Since $q_i$ are $q_j$ coprime, the word $v_i$ is the image of some (special) word  $v_{ij}$ under the morphism $\varphi_{a,q_j}$, plus perhaps some symbols $a$ at the right; and similarly 
$v_j$ is the image of the same word  $v_{ij}$ under $\varphi_{a,q_i}$ (again, with maybe some $a$s at the right).
Thus, $u=\varphi_{a,q_iq_j}(v_{ij})a^r$ for some special word 
$v_{ij}$ with $L(v_{ij})=L(u)$; here $0\leq r < q_iq_j$, and the length of $v$ is
 $\lfloor |u|/q_1q_2 \rfloor$.
 
Extending these arguments to more primes from $Q$, we see that each word $u$ that can be obtained by application of $l$ different morphisms $\varphi_{a,q_{i_1}},\dots,\varphi_{a,q_{i_l}}$ can be obtained by applying all of them consecutively to some word $v$ of length $\lfloor |u|/q_{i_1}
\dots q_{i_l} \rfloor$, and $L(v)=L(u)$.

Thus, by the inclusion-exclusion formula, the contribution to  $s(n)$ of all special words of length
$n\geq 3$ starting with $a$ is
$$
\sum\limits_{l=1}^{k}(-1)^{l+1} \sum\limits_{\{q_{i_1},\dots,q_{i_l}\}\subseteq Q}
s\left( \left \lfloor \frac{n}{q_{i_1}\dots q_{i_l}} \right \rfloor \right ).
$$
The symbol $a$ here is an arbitrary symbol from the alphabet $\Sigma$ of cardinality $d$.
Since $F$ is symmetric with respect to $\Sigma$, we get a recurrent formula 
\begin{equation}
s(n)=d\sum\limits_{l=1}^{k}(-1)^{l+1}
\sum\limits_{\{q_{i_1},\dots,q_{i_l}\}\subseteq Q}
s\left( \left \lfloor \frac{n}{q_{i_1}\dots q_{i_l}} \right \rfloor \right ) 
\quad \text{при~} n \geq 3.
\end{equation}
\begin{example}{\rm
If the alphabet is binary and $Q=\{3,5\}$, then (3) looks like 
$$s(n)=2\left ( s\left ( \left \lfloor \frac{n}{3} \right \rfloor \right )+
s \left ( \left \lfloor \frac{n}{5} \right \rfloor \right )-
s\left ( \left \lfloor \frac{n}{15} \right \rfloor \right ) \right )  $$
and holds for all $n \geq 3$.
}
\end{example}
Note also that for each given $Q$, the values of $s(n)$  for $n=0,1,2,3$ can be found by hand, from the values of $p(n)$ for $n=0,\dots,4$, and depend only on the cardinality of the alphabet and on the fact if 2 and 3 belong to $Q$.

\begin{lemma}
The values of $s(n)$ for $n=0,\dots,3$ are described by the following table.
\begin{center}
\begin{tabular}{|c||c|c|c|c|} \hline
          $n$    & $0$     & $1$        & $2$              & $3$                 \\ \hline \hline

$2 \in Q$,
$3 \in Q$      & $d-1$  & $d^2-d$   & $d^3-d^2$     &$2d^3-3d^2+d$                        \\ \hline

$2 \in Q$,
$3 \notin Q$ &$d-1$   & $d^2-d$  & $d^3-d^2$      & $d^3-d^2$                               \\ \hline

$2 \notin Q$,
$3 \in Q$     & $d-1$   & $d^2-d$  & $2d^2-2d$      & $d^3-d^2$                             \\ \hline

$2 \notin Q$,
$3 \notin Q$ &$d-1$   & $d^2-d$  & $2d^2-2d$      & $d^2-d$                                \\ \hline
\end{tabular}
\end{center}
\end{lemma}

The proof is carried out by direct computation of the function $p(n)$ for  $n=0,\dots,4$.

\section{The growth of $s(n)$}
In this section we use combinatorial techniques to prove several properties of $s(n)$ directly following from (3). For the sake of convenience, let us define the function $s(x)$ for all non-negative reals by the equalities
\begin{equation}
s(x)=s(\lfloor x \rfloor).
\end{equation}
Also, for each finite set $R=\{r_1,\dots,r_m\}$
of primes and for all functions $f(x)$ satisfying (4), define functional operators 
$$
D_R f(x)=\sum\limits_{i=1}^{m}(-1)^{i+1} \sum\limits_{\{s_1,\dots,s_i\}\subseteq R}
f\left( \frac{x}{s_1\dots s_i} \right
)
$$
and
$$
E_R f(x)=f(x)-D_R f(x).
$$
In this terms, (3) can be rewritten as 
\begin{equation}
s(x)=dD_Q s(x) \quad \mbox{~for all~} x \geq 3.
\end{equation}
Note that a solution of this equation satisfying (4) is completely determined by the values of $s(x)$ for $x=0,1,2$.  On the other hand, the solutions of (5) satisfying (4) constitute a linear space: for each pair of such solutions $s_1$ and $s_2$ the function
$\alpha s_1(x)+\beta s_2(x)$ is also a solution and also satisfies (4). 

For each function $f$ we have
$D_{\emptyset}f(x)=0$,
$E_{\emptyset}f(x)=f(x)$ and $D_{\{r\}}f(x)=f(x/r)$.
Moreover, it can be checked directly that
\begin{equation}
E_{R\cup \{r\}}f(x) \equiv E_Rf(x) - E_R f(x/r)
\end{equation}							
for each $r \notin R$, and
\begin{equation}
D_{\{r_1,\dots,r_m\}}f(x)=E_{\emptyset}f(x/r_1)+E_{\{r_1\}}f(x/r_2)+\dots
+ E_{\{r_1,\dots,r_{m-1}\}}f(x/r_m).
\end{equation}
\begin{lemma}
Let $s(x)$ be a solution of {\rm (5)} for the set 
$Q$, and let it satisfy {\rm (4)}. Suppose that 
 $s(3)\geq s(2)\geq s(1)\geq s(0)>0$.
Then the function $s(x)$ is non-decreasing for all $x>0$.
\end{lemma}
{\sc Proof.} Let us define the function $t(x)=s(x)-s(x-1)$.
By the lemma assertion,
$t(x)\geq 0$ for $1\leq x \leq 3$; our goal is to prove that 
$t(n)\geq 0$ for all $n \geq 4$. We shall need two auxiliary statements. 

\begin{proposition}
Let us consider an integer $n>3$ and a subset $Q(n)$ of  $Q$
consisting precisely of the elements of $Q$ dividing  $n$. Then
\begin{equation}
t(n)=dD_{Q(n)}t(n).
\end{equation}
\end{proposition}
{\sc Proof.} As it follows from the definitions, we have
$$
t(n)=s(n)-s(n-1)=
d\sum\limits_{m=1}^{k}(-1)^{m+1}\hskip-6pt \sum\limits_{\{r_1,\dots,r_m\}\subseteq Q}
\left( s\left(\frac{n}{r_1\dots r_m} \right )-s\left(\frac{n-1}{r_1\dots r_m}
 \right )\right ).
$$
Due to (4) the difference 
$(s(n/(r_1\dots r_m))-s((n-1)/(r_1\dots r_m))$
can be non-zero only if $n/(r_1\dots r_m)$ is an integer, that is, if 
$r_1,\dots, r_m \in Q(n)$ (since all the elements of $Q$ are mutually coprime). In this case $s((n)/(r_1\dots r_m))-s((n-1)/(r_1\dots r_m))=t((n)/(r_1\dots r_m))$.

\smallskip
In particular, if some $n>3$ is not multiple to any of the elements of $Q$, then
$t(n)=0$.

\begin{proposition}
Consider an integer $n>3$, $n=q_1^{i_1}\dots q_l^{i_l}n_0$, where
$\{q_1,\dots, q_l\} = Q(n)$  and $n_0$ is coprime with all elements of 
$Q$. Suppose that $t(n_0)\geq 0$; then $t(n)\geq 0$.
\end{proposition}
{\sc Proof.} Let us proof a bit more, namely, that 
$E_{\{q_1,\dots,q_j\}}t(n)\geq 0$ for all $j=0,\dots,l$ (in particular, for 
$E_{\emptyset}t(n)=t(n)$). Let us use induction on
$i=i_1+\dots +i_l$; the base of induction for $i=0$, that is, for $n=n_0$,
is asserted since $Q(n_0)=\emptyset$.

The induction step is based on combining the equalities (8) and (7), so that
$$
t(n)=E_{\emptyset}t(n)=dE_{\emptyset}t(n/q_1)+\dots +d
E_{\{q_1,\dots,q_{l-1}\}}t(n/q_l)
$$
and
\begin{multline*}
E_{\{q_1,\dots,q_j\}}t(n)=(d-1)E_{\emptyset}t(n/q_1)+\dots
+(d-1)E_{\{q_1,\dots,q_{j-1}\}}t(n/q_j)
\\
+
dE_{\{q_1,\dots,q_{j}\}}t(n/q_{j+1})+\dots+
dE_{\{q_1,\dots,q_{l-1}\}}t(n/q_l)
\end{multline*}
for all $j=1,\dots,l$. By the definitions, for all $m$ we see that $n/q_m$ is an integer, the set 
$Q(n/q_m)$ is equal either to $Q(n)$ or to $Q(n)\backslash \{q_m\}$,
and the sum of degrees of 
$q$s 
from $Q(n/q_m)$ 
is $i-1$. So, by the induction hypothesis, the proposition can already be applied for all the operators 
$E_R$ 
from the right part of the system above, and the value of each  
$E_{\{q_1,\dots,q_j\}}t(n)$
 is a sum of non-negative summands. In particular, this is true for $j=0$, that is, for the function $t(n)$ itself. 
 
\smallskip
To prove the lemma, it is now sufficient to note that Proposition 2 can be applied to each $n>3$. Indeed, $n_0$
is equal either to $1,2$ or $3$ (and then $t(n_0)\geq 0$ by the assertion),
or it is a number greater than 3, coprime to all elements of $Q$  (then $t(n_0)=0$ by Proposition 1).

\begin{lemma}
Let the function $s(x)$ be a solution of {\rm (5)} satisfying
{\rm (4)},
and suppose that $s(x)>0$ for all $x\geq 0$ and $s(0)<s(1)<s(2)$.
Then there exists a positive constant $C$ such that for all $x>1$ the inequality holds
$$s(x) \leq Cx^{\alpha},$$
where $\alpha$ is the solution of  {\rm (1)}.
\end{lemma}
{\sc Proof.} In this proof, it is convenient to assume that the elements $Q$ are listed in ascending order.

Note that the function $s(x)$ and the operators  
$E_{\{q_1,\dots,q_i\}}s(x)$ for $1\leq i \leq k$ take only a finite number of values when
$1 \leq x \leq q_1\dots q_k$, so that there exists a positive constant 
$C$ such that 
\begin{equation}
\begin{split}
s(x)&\leq  C x^{\alpha}\\
E_{\{q_1\}}s(x)&\leq  (1-1/q_1^\alpha)C x^\alpha\\
& \vdots \\
E_{\{q_1,\ldots,q_{k-1}\}}s(x)& \leq  (1-1/q_1^\alpha)\cdots (1-1/q_{k-1}^\alpha)C x^\alpha
\end{split}
\end{equation}
for all $1 \leq x < q_1\dots q_k$.
This system of inequalities constitute the base of induction; for the induction step, suppose that (9) holds for all 
$x \in [1,X)$, where $X \geq q_1\dots q_k$, and prove it for $x \in [X,q_1X)$ (note that by our assertion, $q_1$ is the minimal element of $Q$).
Due to the definitions of $D_Q$ and $E_Q$, and also to 
(5)--(7), for all $x\geq 3$ we have

\begin{multline}
E_{\{q_1,\dots,q_{k-1}\}}s(x)=(d-1)
[s(x/q_1)+E_{\{q_1\}}s(x/q_2)
+\dots + E_{\{q_1,\dots,q_{k-2}\}}s(x/q_{k-1}) ]
\\
+ dE_{\{q_1,\dots,q_{k-1}\}}s(x/q_k).
\end{multline}
Substituting the inequalities from (9), which hold by the induction hypothesis, to the operators in the right part of this equality, we obtain that
\begin{multline*}
E_{\{q_1,\dots,q_{k-1}\}}s(x)\leq Cx^{\alpha}(d-1)
\left [{1\over q_1^{\alpha}}+
\biggl(1-{1\over q_1^{\alpha}}\biggr)
{1\over q_2^{\alpha}}+\dots+
{1\over q_{k-1}^\alpha}
\prod\limits_{i=1}^{k-2}\biggl(1-{1\over q_i^\alpha}\biggr) \right ]
\\
+  Cx^{\alpha}d \cdot 
{1\over q_{k}^\alpha}\prod\limits_{i=1}^{k-1}\biggl(1-{1\over q_i^\alpha}\biggr)
\\
= Cx^{\alpha}(d-1)
\sum\limits_{j=1}^{k}
{1\over q_{j}^\alpha}\prod\limits_{i=1}^{j-1}
\biggl(1-{1\over q_i^\alpha}\biggr) 
+ Cx^{\alpha} \cdot {1\over q_{k}^\alpha}
\prod\limits_{i=1}^{k-1}\biggl(1-{1\over q_i^\alpha}\biggr).
\end{multline*}
Note that 
$$
\sum\limits_{j=1}^{k}{1\over q_{j}^\alpha}\prod\limits_{i=1}^{j-1}
\biggl(1-{1\over q_i^\alpha}\biggr)
=1-\prod\limits_{i=1}^{k}\biggl(1-{1\over q_i^\alpha}\biggr).
$$
So, using (1), we see that the first summand in the right part of the previous inequality is 
$$
Cx^{\alpha}(d-1)[1-(d-1)/d]=Cx^{\alpha}(d-1)/d.
$$
In its turn, the second summand can be transformed as follows: 
\begin{multline*}
Cx^{\alpha} \cdot {1\over q_{k}^\alpha}
\prod\limits_{i=1}^{k-1}
\biggl(1-{1\over q_i^\alpha}\biggr)
=Cx^{\alpha} 
\Biggl[\prod\limits_{i=1}^{k-1}
\biggl(1-{1\over q_i^\alpha}\biggr)
-\prod\limits_{i=1}^{k}\biggl(1-{1\over q_i^\alpha}\biggr)\Biggr]
\\
=Cx^{\alpha} 
\Biggl[\prod\limits_{i=1}^{k-1}\biggl(1-{1\over q_i^\alpha}\biggr)-
{d-1\over d}\Biggr].
\end{multline*}
The sum of these two summands gives precisely the right part of the last inequality of the system
(9), which is what we needed. Thus, the last inequality from  (9)
holds also for $x\in [X,q_1X)$.

Now to prove all the previous inequalities from (9) for $x \in [X,q_1X)$, it is sufficient to use consecutively in descending order, for $i=k-2,\dots,0$, the inequalities 
\begin{equation}
E_{\{q_1,\dots,q_i\}}s(x)=E_{\{q_1,\dots,q_{i+1}\}}s(x)+E_{\{q_1,\dots,q_i\}}s(x/q_{i+1}),
\end{equation}
which hold due to (6). Substituting to them respective inequalities from (9), valid by the induction hypothesis, we see that the inequalities from (9), and in particular the first of them, hold also for $x\in [X,q_1X)$.

\section{Lemma of uniqueness of the root}
\begin{lemma}
For each finite set $Q$ of primes of cardinality at least 2 the equation
$$
\frac{d-1}{d}=\prod\limits_{q\in Q} ( 1-q^{-x} )
$$
have no roots with the real part equal to $\alpha=\alpha(Q,d)$ except for the root $x=\alpha$.
\end{lemma}
{\sc Proof.} Suppose that there exists another root $\alpha+iy$, $0\ne y \in {\mathbb  R}$. Then
$$
\frac{d-1}{d}=\prod\limits_{q \in Q} ( 1-q^{-\alpha-iy} )
=\prod\limits_{q \in Q} |1-q^{-\alpha-iy} |.
$$
By the triangle inequality we have $|1-q^{-\alpha-iy}|\geq 1-q^{-\alpha}$,
and the equality is reached if and only if $q^{-iy}=1$,
that is, $y=2\pi k/\ln q$, $k \in \mathbb  Z$.
This equality is possible for at most one prime  $q$. Indeed, if 
$y=2\pi k_1/\ln q_1=2\pi k_2 / \ln q_2$, then
$q_2^{k_1} = q_1^{k_2}$, which is impossible.
Multiplying the inequalities $|1-q^{-\alpha-iy}|\geq 1-q^{-\alpha}$
for $q \in Q$ (at least one of them is strict since $Q$ contains at least two primes), we get 
$$
\prod\limits_{q \in Q} |1-q^{-\alpha-iy} |>
\prod\limits_{q \in Q} |1-q^{-\alpha} |=\frac{d-1}{d};
$$
a contradiction. 

\medskip
Note that it is important here that $Q$ contains at least two elements. For $|Q|=1$, the lemma does not hold.

\section{Tauberian theorem}

Denote by $W$ the Wiener class of functions $g$ summable on $(0,+\infty)$ and such that the Mellin transform 
 $\int\limits_0^{\infty} g(t)t^{ix}\,dt$ is non-zero for all real $x$.

\begin{theorem}[Wiener-Pitt, \cite{hardi}, Theorem 233]
Suppose that $g\in W$, and the function $f$ is real-valued, bounded and slowly decreases in the sense that if $y>x\rightarrow +\infty$, $y/x\rightarrow 1$, then
$\liminf f(y)-f(x)\geq 0$. Then the fact that
$$
\frac1x \int\limits_0^{\infty} f(t)g(t/x) dt\rightarrow l\int\limits_0^{\infty} g(t)\,dt
$$
for $x\rightarrow +\infty$ implies that $f(t)\rightarrow l$
for $t\rightarrow +\infty$.
\end{theorem}
\begin{theorem}
Suppose that real numbers $c_1$, $c_2$, $\dots$, $c_l$
and $m_1$, $m_2$, $\dots$, $m_l$ $(m_i>1$ for all $i)$ are such that the equation 
$F(\alpha):=\sum c_i m_i^{-\alpha}=1$
has a unique positive solution $\alpha=\alpha_0$, and $F(\alpha_0+ix)\neq 1$ for any real $x\ne 0$ and
$F'(\alpha_0)\ne 0$. Suppose also that a non-decreasing function $s(x)$, $x\geq 0$,
such that $s(x)=s(\lfloor x \rfloor)$, satisfies the equation
$$
s(x)=\sum\limits_{i=1}^{l} c_i s(x/m_i), \quad x\geq N_0>0,
$$
and a relation $s(x)=O(x^{\alpha_0})$ for $x\rightarrow +\infty$.

The the quotient $s(x)/x^{\alpha_0}$ tends to some limit with $x \to \infty$.
\end{theorem}
{\sc Proof.} Denote by $\chi(t)=\chi_{[0,1]}(t)$ the characteristic function of the interval $[0,1]$ and consider the kernel
$$
g(t)=\frac{\chi(t)-\sum c_im_i^{-\alpha_0} \chi(m_it)}t, \quad
t>0.
$$
Note that $g(t)$ is a piecewise constant finite function whose support is separated from zero and infinity. So, 
 $g\in L_1(0,+\infty)$.
Let us calculate its Mellin transform:
$$
\int\limits_0^{\infty} g(t)t^{ix} dt=\frac{1-\sum c_im_i^{-\alpha_0-ix}}{ix};
$$
the value for $x\ne 0$ is defined by continuity. Our assertion on the properties of $F(\alpha)$ can now be stated as follows:
the Mellin transform of $g$ has no real zeros. Thus, the function  $g$ belongs to the Wiener class $W$ on $(0,+\infty)$.

Note that the expression
\begin{equation}
I(x)=\frac1x \int\limits_0^{\infty} \frac{s(t)}{t^{\alpha_0}} g(t/x) \,dt
\end{equation}
does not depend on $x$ when $x$ is positive and sufficiently large, and is then equal to 
\begin{equation}
A=\sum c_i m_i^{-\alpha_0} \int\limits_{N_0/m_i}^{N_0}
\frac{s(\tau)}{\tau^{\alpha_0+1}} \,d\tau.
 \end{equation}

Indeed, the integrand in (12) is equal to 0 for $t<N_0$ (for $x$ sufficiently large), so that the lower limit of integration can be replaced by  $N_0$. So, the integral can be rewritten as 
$$
I(x)=\int\limits_{N_0}^x \frac{s(t)}{t^{\alpha_0+1}} \,dt
-\sum c_i m_i^{-\alpha_0} \int\limits_{N_0}^{x/m_i}
\frac{s(t)}{t^{\alpha_0+1}} \,dt.
$$
Let us represent the first integral as the sum 
$s(t)=\sum c_i s(t/m_i)$ for $t\geq N_0$
and replace in the respective summands $t$ to $\tau=t/m_i$. We obtain (13).

Note that if $s(x)$ is increasing and  $s(x)=O(x^{\alpha_0})$,
then the bounded function $f(x)=s(x)\cdot x^{-\alpha_0}$
is also slowly decreasing. Indeed,
$$
f(y)-f(x)=s(y)y^{-\alpha_0}-s(x)x^{-\alpha_0}
\geq s(x)(y^{-\alpha_0}-x^{-\alpha_0})=
(s(x)x^{-\alpha_0})((x/y)^{\alpha_0}-1),
$$
where the first multiple is bounded and the second tends to zero. 

So, by the Tauberian Wiener-Pitt theorem,
$s(x)/x^{\alpha}=f(x)\rightarrow A/\int g
=A/(-F'(\alpha_0))$ for $x \to \infty$,
which was to be proved.

\section{End of the proof}

Note that sometimes Theorem 3 cannot be applied directly to $s(x)$ since this function is not non-decreasing, as we see from Lemma 2. However, as it was mentioned in the beginning of Section 5, the solutions of (5) satisfying (4) constitute a linear space, so that  $s(x)$ can be represented as a difference  $s(x)=s_1(x)-s_2(x)$, where $s_1$ are $s_2$ also solutions of (5) satisfying (4), and they are non-decreasing for $x\leq 4$. Due to Lemma 3, they are non-decreasing for all positive numbers, and we can apply Theorem 3 (and all the previous statements necessary for it) to each of them separately.

The coefficients $c_r$ and $m_r$ from the statement of Theorem 3 are defined naturally by (3): for each $m_r$ equal to  the product $q_{i_1}\dots q_{i_l}$, the respective $c_r$ is $d(-1)^{l+1}$. It is not difficult to check that (1) can be rewritten in this notation as it is needed for the theorem; here the function $F(x)$ is defined by 
$$
F(x)=d\left ( 1- \prod\limits_{i=1}^{k} \bigl( 1-q_i^{-x} \bigr) \right ).
$$
It can be easily seen that its derivative at the point $x=\alpha$ is negative; the remaining conditions of the theorem are fulfilled due to Lemmas 4 and 5. Thus, due to Theorem 3, for $x \to \infty$ there exists the limit of the ratio
$s(x)/x^{\alpha}$: 
or, strictly speaking, the limits of $s_1(x)/x^{\alpha}$ and $s_2(x)/x^{\alpha}$, but then we can just take their difference to get the needed limit of $s(x)/x^{\alpha}$.

It remains to note that $s(n)$ is the first difference function of the arithmetical complexity $a_w(n)=p_{Q,d}(n)$ we are really interested in. For it, the limit $p_{Q,d}(n)/x^{\alpha+1}$ exists by the Stolz–Ces\`aro theorem [14] and is equal to the limit for the first differences divided by 
$\alpha+1$. Let us write it down explicitly.

First let us compute the constant $A$ from (13): it can always be done directly for 
 $s$ not for $s_1$ and $s_2$, since the integral of a linear combination is the linear combination of integrals. This constant depends on the fact if 2 belongs to $Q$. Suppose first that 
$2 \not\in Q$. In this case $\sum c_i=d$. The values of $N_0/m_i$
are not greater than 1 for all $i$, so that 
\begin{eqnarray*}
A&=&\alpha^{-1}\sum c_i m_i^{-\alpha}\left(s(0)(3^{-\alpha}m_i^{\alpha}-1)+s(1)(1-2^{-\alpha})+s(2)(2^{-\alpha}-3^{-\alpha})\right)= \\
&&\alpha^{-1}\left(3^{-\alpha}(s(0)d-s(2))+2^{-\alpha}(s(2)-s(1))+s(1)-s(0) \right).
\end{eqnarray*}
It remains to substitute into this expression the values $s(0)=d-1$, $s(1)=d^2-d$, $s(2)=2d^2-2d$ to get
$$
A=\frac1{\alpha}((d-1)^2+(2^{-\alpha}-3^{-\alpha})(d^2-d)).
$$

Analogously, for $2 \in Q$ we obtain that 
$$
A=\frac{1}{\alpha}(d-1)^2.
$$

So, the constant equal to the limit of 
$p_{Q,d}(n)/x^{\alpha+1}$ can be found directly and is equal to $\frac{A}{-F'(\alpha)(\alpha+1)}$. In particular we see that it is positive, which completes the proof of the theorem. 

\section{Acknowledgements}

We are greateful to professors V. S. Guba and F. L. Nazarov for the useful discussions and to the maintainers of the LiveJournal community {\tt ru-math} giving a rare opportunity of interdisciplinary contacts of mathematicians.


\begin{thebibliography}{10}
\bibitem{ferenczi}
Ferenczi S.,
Complexity of sequences and dynamical systems //
Discrete Math. 206, 145--154, 1999.
\bibitem{cas_int}
Cassaigne J.,
Constructing infinite words of intermediate complexity //
LNCS 2450, Springer Verlag, 173--184, 2002.
\bibitem{aff}
Avgustinovich S. V., Fon-Der-Flaass D. G., Frid A. E., Arithmetical complexity of infinite words. // Masami Ito and Teruo Imaoka, editors, Words, Languages and Combinatorics III, pp. 51-62, Singapore, 2003.
\bibitem{symm}
Frid A. E.,   Arithmetical complexity of symmetric D0L words //  Theoret. Comput. Sci. 306 (2003) 535--542.
\bibitem{lin}
Frid A. E., Sequences of linear arithmetical complexity //  Theoret. Comput. Sci. 339 (2005) 68--87.
\bibitem{nlogn}
Frid A. E., On possible growths of arithmetical complexity // Theoret. Informatics
Appl., 2006, 40, no. 3, 443--458.
\bibitem{salimov}
Salimov P. V., Constructing infinite words of intermediate arithmetical complexity // LNCS 5457, Springer Verlag, 696--701, 2009.
\bibitem{bes}
Besicovitch A. S., On the linear independence of fractional powers of integers // J. London Math. Soc. 15 (1940), 3--6.
\bibitem{schanuel}
Waldschmidt M., Open Diophantine problems // Mosc. Math. J. 4 (2004), 245--305.
\bibitem{koskas}
Koskas M., Complexit\'es de suites de Toeplitz // Discrete Math. 183 (1998)
161--183.
\bibitem{ck}
Cassaigne J., Karhum\"aki J., Toeplitz words, generalized periodicity and
periodically iterated morphisms // European J. of Combinatorics 18 (1997),
497--510.
\bibitem{cas!}
Cassaigne J., Complexit\'e et facteurs sp\'eciaux // Bull. Belg. Math. Soc. Simon Stevin 4 (1997),
67--88.
\bibitem{hardi}
G. Hardy, Divergent series. ...
\bibitem{fiht}
G. M. Fichtenholz, ???







\end{thebibliography}
\end{document}